\documentclass[journal]{IEEEtran}
\pdfoutput=1
\usepackage{cite}

\usepackage{float}

\ifCLASSINFOpdf
\else
\fi
\usepackage{amssymb}
\usepackage{amsmath}

\usepackage{graphicx}
\begin{document}

%
\title{Observations towards a proof of the Collatz conjecture using $2^{j}k+x$ number series}
%
%
%

\author{J. Stöckl
}

%
%

\markboth{arXiv notes}%
{}
%



\maketitle

\begin{abstract}

The document tries to put focus on sequences with certain properties and periods leading to the first value smaller than the starting value in the Collatz problem. With the idea that, if all starting numbers lead ultimately to a smaller number, all full sequences lead to 1 with a finite stopping time, the problem could be reduced to more structured shorter sequences. It is shown that this sequences exist and follow consistent rules. Potential features of an infinite cycle, also leading to a smaller number, are also discussed. Further, an argument for only one possible closed cycle is given for the special sequence of alternating odd and even steps as well as arguments that infinite cycles must exist. Using the observation that periodic behavior is exists an additional argument is provided the probability of subsets which will end up at a number smaller the initial value possibly even of $\mathbb{N}^{*}$ will indeed end up at unity in the Collatz problem \cite{MANCOSU_2009}. The work is to be seen as in progress and shared as an contribution to discussion rather than a concrete publication.

\end{abstract}


%
\IEEEpeerreviewmaketitle

\section{Introduction}
%
%
%
%
\IEEEPARstart{T}{he} Collatz conjecture is one of the best known unsolved mathematical problems and is therefore only briefly described as follows.
The sets of natural numbers are given by:

\[\mathbb{N}:= \left \{ 0,1,2,3... \right \}\]

\[\mathbb{N}^{*}:= \left \{ 1,2,3... \right \}\]

Be \(n\in \mathbb{N}^{*}\) and \(Col:\mathbb{N}^{*}\rightarrow \mathbb{N}^{*}\) the Collatz function then

\[Col(n)=\left\{\begin{matrix}
\frac{n}{2} \; if \; n\equiv 0\;(mod \; 2) \;later\; named \;E\\ 
3n+1 \; if \; n \equiv 0 \;(mod \; 2)\;later\; named \;O\
\end{matrix}\right.\]

\begin{figure}[H]
    \centering
    \includegraphics[width=0.5\linewidth]{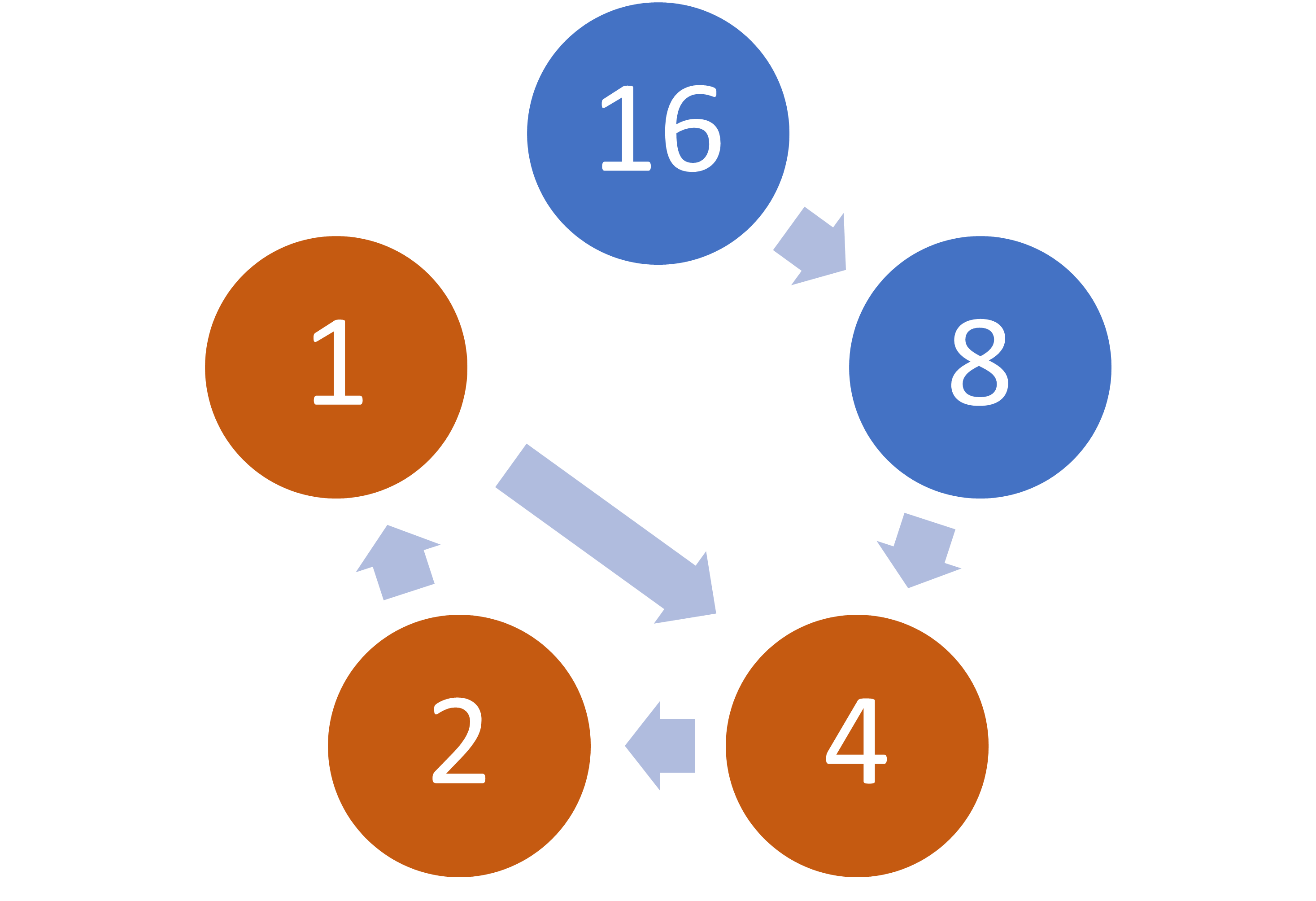}
    \caption{Basic cycle ending with 1-4-2-1}
    \label{fig:basic_cycle}
\end{figure}

The Collatz conjecture states that a sequence of application of the Collatz function always ends with the number 1 or more specifically with the sequence 1, 4, 2, 1, etc. Exemplarily, this is shown with the sequence starting from 16 (Fig. \ref{fig:basic_cycle}). Experimental evidence using numerical check with computers have been provided up to \(2^{68}\approx 2.95\times10^{20}\) \cite{barina_convergence_2021}. Previous research has also put emphasis on the identification of simple cycles such as Mimuro \cite{mimuro}. Here, the focus is put on the properties of such cycles with respect to the whole available set of \(\mathbb{N}^{*}\) as starting numbers. In the context of diophantine equations work has been published by e.g. Belaga \textit{et al.} \cite{belaga}. Also mapping to number systems has been investigated \cite{IDOWU2015105}.

\begin{figure}
    \centering
    \includegraphics[width=0.9\linewidth]{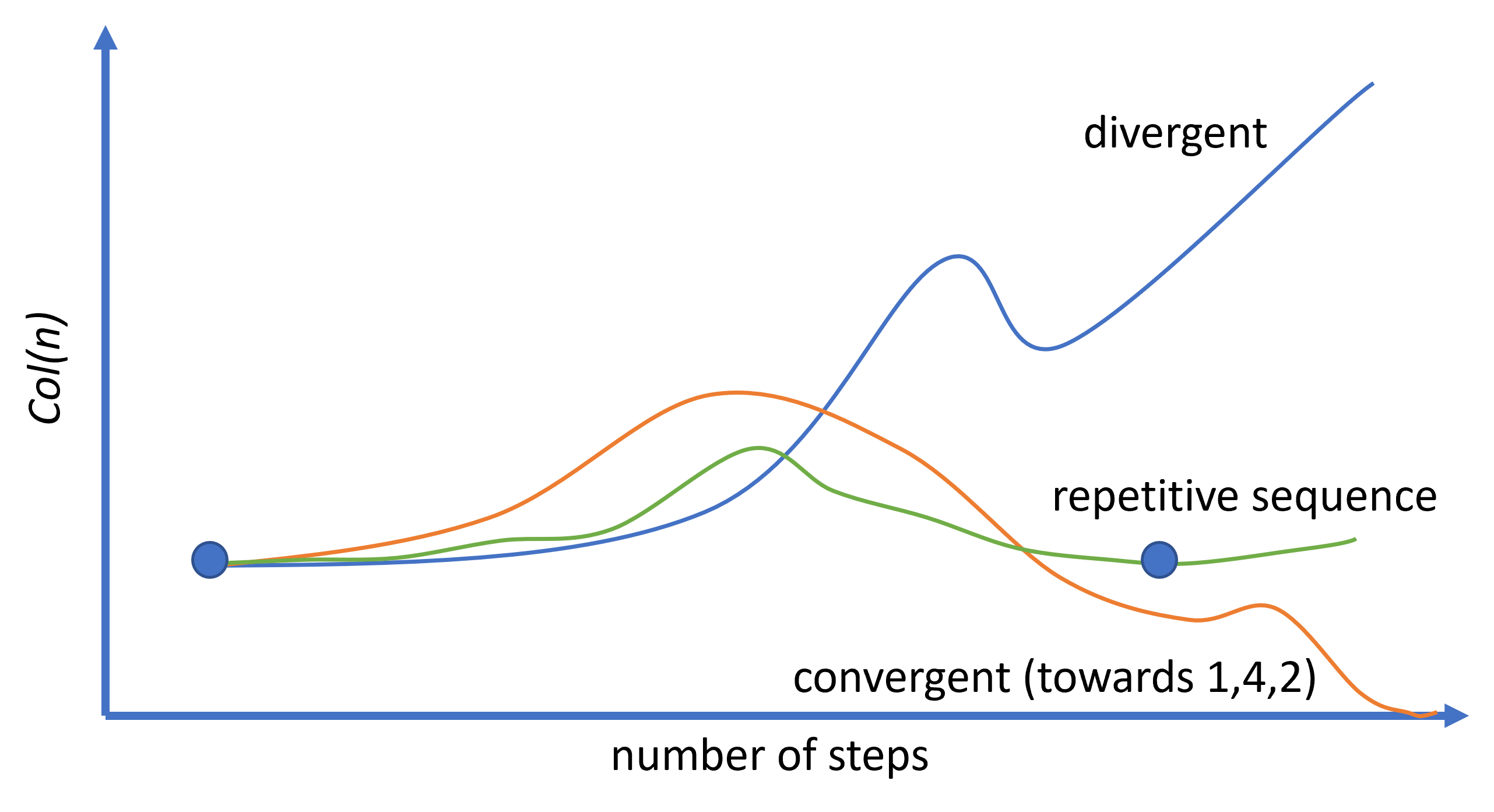}
    \caption{Possible behaviours of Collatz sequences}
    \label{fig:schematic-behaviour}
\end{figure}

All sequences will have one of the properties shown in Fig. \ref{fig:schematic-behaviour}, although only the convergent case has been observed in the Collatz problem.

Although in general the stopping time, means the number of steps needed to reach the number 1, is the focus, here the number of steps reaching a number smaller (or equal) than the starting number is of interest. The question, which is in principle the Collatz conjecture rephrased, is whether all numbers \(n\in \mathbb{N}^{*}\) (exept 1) will reach this lower number after a certain "lower number time". This is insofar important as it will be shown that there are patterns to be identified which would not evolve further and enable a new way of clustering and analysis of the numbers in the context of the Collatz conjecture. The stopping time is then simply the sum of all individual lower number times.
In the following starting from the most trivial cases a set of rules will be developed to make further assumptions and some basic proofs. Ultimately, it should be possible to justify the statements for all natural numbers.

The notes start with basic cycles for odd and specific even numbers to make the point. Then the observations are extended to the next longer cycles of length 6 and 8. General rules and formulas are derived and finally extended to more general cases. With the starting number 27 the case of long cycles including the next possible starting value with this sequence is presented. Finally, potential further investigations are proposed.

\section{Basic Cycles}

\subsection{Basic case: even natural numbers}
The most trivial case are the even natural numbers as already the subsequent number is lower than the starting number (Table \ref{tab:basic_even}).

\[n_{0}\in \begin{Bmatrix}
2k
\end{Bmatrix};k\in \mathbb{N}^{*}\rightarrow Col(n)=n/2 =k< n\]

\begin{table}
    \centering
    \begin{tabular}{|c|c|c|c|c|c|}
    \hline
        Number 1 & Number 2 & Number 3 & General & Step & Cycle \\
        \hline
        4 & 200 & 750 & \textit{n} & 1 & \textit{E}\\
        2 & 100 & 375 & \textit{n/2} &  & \\
        \hline
    \end{tabular}
    \caption{Basic cycle: even numbers}
    \label{tab:basic_even}
\end{table}

This is an important result as this means that the above assumption is true for all even natural numbers and the only pending set is the odd natural numbers without 1:

\[n_{0}\in \begin{Bmatrix}
2k+1
\end{Bmatrix};k\in \mathbb{N}^{*}\]

The following statements can be considered as true:
\begin{itemize}
    \item The last step in a sequence needs to be a division to reach this lowest number and hence, the second to last number needs to be an even natural number.
    \item The starting number needs to be an odd natural number to avoid the trivial case.
\end{itemize}

\subsection{Basic case: even natural numbers}

Following a multiplication with 3 and subsequent division by 2 is still larger than the starting number. This means that at least a second division is needed to reach the goal (Table \ref{tab:basic_odd}).

\begin{table}
    \centering
    \begin{tabular}{|c|c|c|c|c|c|}
    \hline
      Number 1   & Number 2 & Number 3 & General & Step & Cycle \\
         \hline
       5  & 77 & 853 & \textit{n} & 1 & \textit{O}\\
       16  & 232 & 2560 & (3\textit{n}+1) & 2 & \textit{E}\\
        8 & 116 & 1280 & (3\textit{n}+1)/2 & 3 & \textit{E}\\
        4 & 58 & 640 & (3\textit{n}+1)/4 &  & \\
         \hline
    \end{tabular}
    \caption{Basic cylce: odd numbers}
    \label{tab:basic_odd}
\end{table}

It is already known that all starting numbers share this behavior when they follow the construction rule: 
\[n_{0}\in \begin{Bmatrix}
4k+1
\end{Bmatrix};k\in \mathbb{N}\]

This property again halves the remaining set with the following subset remaining:

\[n_{0}\in \begin{Bmatrix}
4k+3
\end{Bmatrix};k\in \mathbb{N}^{*}\]

\section{Sequence Properties}

Which cycles are possible to reach the first smaller number. Obviously too many \textit{O} operations will lead to a much too large number, although an \textit{O} operation is always followed by an \textit{E} operation. Too many \textit{E} operations will lead to a smaller number than necessary going beyond the stopping rule.

Be \textit{i} the number of \textit{O} operations and \textit{j} the number of \textit{E} operations the following rule should be sufficient to identify the number steps and operations needed to be able to succeed.

\[2^{j}-3^{i}> 0\]

or
\[j > i \cdot log_2(3)\]

as $log_2(3)$ is irrational and $i, j$ are natural numbers there can never be unity.

Table \ref{tab:cycle_length} shows the results for up to a cycle length of 37 steps. It can be seen that not all cycle lengths are possible in order to obtain a smaller number.

\begin{table}
    \centering
\begin{tabular}{|c|c|c|c|c|}
\hline
E ops & O ops & Result  & Cycle length & Remark         \\
\hline
1            & 0            & 1       & 1            & Even number    \\
2            & 1            & 1       & 3            & Shortest cycle \\
3            & 2            & -1      & 5            & Not possible   \\
4            & 2            & 7       & 6            & Possible       \\
4            & 3            & -11     & 7            & Not possible   \\
5            & 3            & 5       & 8            & Possible       \\
5            & 4            & -49     & 9            & Not possible   \\
6            & 4            & -17     & 10           & Not possible   \\
7            & 4            & 47      & 11           & Possible       \\
8            & 5            & 13      & 13           & Possible       \\
9            & 6            & -217    & 15           & Not possible   \\
10           & 6            & 295     & 16           & Possible       \\
11           & 7            & -139    & 18           & Not possible   \\
12           & 7            & 1909    & 19           & Possible       \\
13           & 8            & 1631    & 21           & Possible       \\
14           & 9            & -3299   & 23           & Not possible   \\
15           & 9            & 13085   & 24           & Possible       \\
16           & 10           & 6487    & 26           & Possible       \\
17           & 11           & -46075  & 28           & Not possible   \\
18           & 11           & 84997   & 29           & Possible       \\
19           & 12           & -7153   & 31           & Not possible   \\
20           & 12           & 517135  & 32           & Possible       \\
21           & 13           & 502829  & 34           & Possible       \\
22           & 14           & -588665 & 36           & Not possible   \\
23           & 14           & 3605639 & 37           & Possible  \\
\hline
\end{tabular}
    \caption{Cycle length}
    \label{tab:cycle_length}
\end{table}

\section{Further cycles}

\subsection{The cycles for the cycle length 6 and first results}

The next cycle length is 6 and consists of 4 even operations and 2 odd operations. Taking into account the boundary conditions for cycles it is clear that only steps 3, 4, and 5 can be selected and that one \textit{O} and two \textit{E} need to be distributed. This means that three cycles can be defined (see Table \ref{tab:cycle_6}). It can also be seen that cycles 2 and 3 are starting with the sequence of the shortest cycle leading already to a number smaller than the starting number. This means that only cycle 1 can be used to find a smaller number.

\begin{table}
    \centering
    \begin{tabular}{|c|c|c|c|c|}
    \hline
        Cycle 1 & Cycle 2 & Cycle 3  & Step & Remark\\
        \hline
         & shortest cycle & shortest cycle &  & \\
        \textit{O} & \textit{O} & \textit{O} & 1 & needs to be \textit{O}\\
       \textit{E}  & \textit{E} & \textit{E} & 2 & needs to be \textit{E}\\
        \textit{O} & \textit{E} & \textit{E} & 3 & free\\
        \textit{E} & \textit{O} & \textit{E} & 4 & free\\
        \textit{E} & \textit{E} & \textit{O} & 5 & free\\
        \textit{E} & \textit{E} & \textit{E} & 6 & needs to be \textit{E}\\
        possible & not possible & not possible &  & \\
        \hline
    \end{tabular}
    \caption{Sequences with length 6}
    \label{tab:cycle_6}
\end{table}

Table \ref{tab:length6_numbers} shows a set of example numbers with a first generalization formula.

\begin{table}
    \centering
    \begin{tabular}{|c|c|c|c|c|c|}
    \hline
        N 1 & N 2 & N 3 & General & Step & Cycle\\
         \hline
        3 & 19 & 163 & \textit{n} & 1 & \textit{O}\\
        10 & 58 & 490 & (3\textit{n}+1) & 2 & \textit{E}\\
        5 & 29 & 245 & (3\textit{n}+1)/2 & 3 & \textit{O}\\
        16 & 88 & 736 & 3(3\textit{n}+1)/2+1 & 4 & \textit{E}\\
        8 & 44 & 368 & (3(3\textit{n}+1)/2+1)/2 & 5 & \textit{E}\\
        4 & 22 & 184 &(3(3\textit{n}+1)/2+1)/4 & 6 & \textit{E}\\
        2 & 11 & 92 & (3(3\textit{n}+1)/2+1)/8 &  & \\
         \hline
    \end{tabular}
    \caption{Selected sequences with length 6}
    \label{tab:length6_numbers}
\end{table}

Again, the numbers share a property which can be identified with the general formula namely:

\[n_{0}\in \begin{Bmatrix}
16k+3
\end{Bmatrix};k\in \mathbb{N}\]

It is now time to make a next observation with respect to the numbers and their cycles (see Table \ref{tab:sequence_subsets}). All sequences so far share a common structure based on the power of 2 plus an offset value.

\begin{table}
    \centering
    \begin{tabular}{|c|c|c|c|}
    \hline
        Cycle length & Subset & \textit{E} cycles & General subset\\
         \hline
        1 & 2,4,6... & 1 & \(2^{1}k\)\\
        3 & 1,5,9... &  2& \(2^{2}k+1\)\\
        6 & 3,19,35... & 4 & \(2^{4}k+3\)\\
         \hline
    \end{tabular}
    \caption{First sequence subsets}
    \label{tab:sequence_subsets}
\end{table}

\subsection{The cycles for the cycle length 8}

With this first result, we clearly hope for something similar to Table \ref{tab:length8_anticipation}

\begin{table}
    \centering
    \begin{tabular}{|c|c|c|c|}
    \hline
        Cycle length & Subset & \textit{E} cycles & General subset\\
        \hline
        8 & \textit{x}, 32+\textit{x}, 64+\textit{x}, 96+\textit{x} & 5 & \(2^{5}k+x\) \\
        \hline
    \end{tabular}
    \caption{Length 8 anticipation}
    \label{tab:length8_anticipation}
\end{table}

It turns out that, although more cycle combinations would be possible, there are previous cycles as subsets leading already to an earlier lower number. Ultimately, only two cycles remain as possible options available as seen in Table \ref{tab:sequence8_length}.

\begin{table}
    \centering
    \begin{tabular}{|c|c|c|c|}
    \hline
       Cycle 1  & Cycle 2 & Step & Remark\\
         \hline
        \textit{O} & \textit{O} & 1 & needs to be \textit{O} \\
        \textit{E} & \textit{E} &  2& needs to be \textit{E}\\
        \textit{O} & \textit{O} &  3& needs to be \textit{O}\\
        \textit{E} & \textit{E} &  4& depends\\
        \textit{E} & \textit{O} &  5& depends\\
       \textit{O}  & \textit{E} &  6& depends\\
        \textit{E} & \textit{E} &  7& needs to be \textit{E}\\
        \textit{E} & \textit{E} &  8& needs to be \textit{E}\\
       possible  & possible &  & \\
         \hline
    \end{tabular}
    \caption{Sequences with length 8}
    \label{tab:sequence8_length}
\end{table}

Examples for the two cycles can be found in the Appendix. As a result it can be seen in Table \ref{tab:general_length8} that the anticipated behavior is fulfilled.

\begin{table}
    \centering
    \begin{tabular}{|c|c|c|c|}
    \hline
       Cycle length  & Subset & \textit{E} cycles & General\\
         \hline
       8  & 11,32+11,64+11,96+11,... & 5 & \(2^{5}k+11\)\\
        8 & 23,32+23,64+23,96+23,... & 5 & \(2^{5}k+23\)\\
         \hline
    \end{tabular}
    \caption{The starting values for the sequences with the length of 8 steps}
    \label{tab:general_length8}
\end{table}

\begin{figure}
    \centering
    \includegraphics[width=0.9\linewidth]{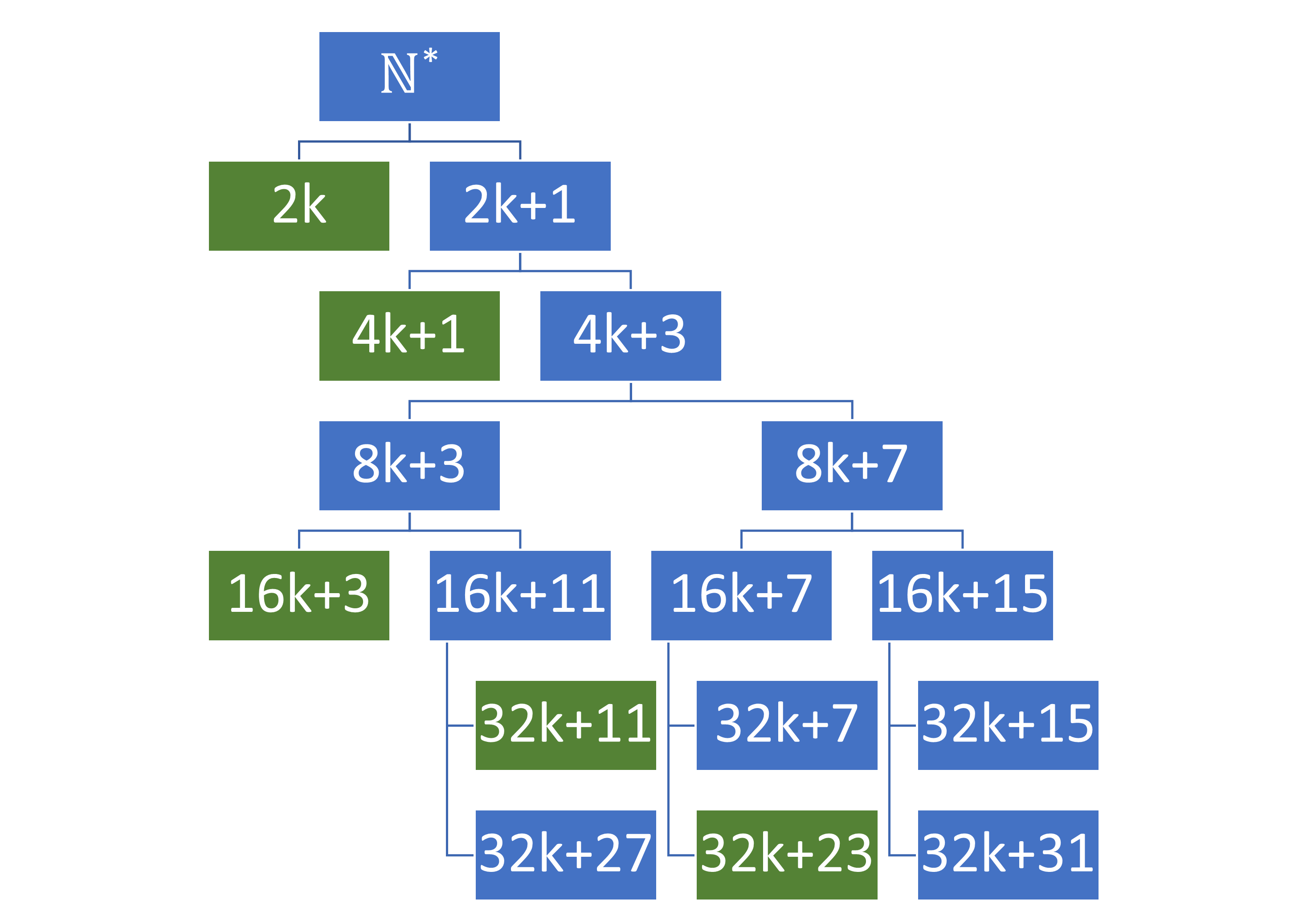}
    \caption{Subsets with sequence (green) and without sequence (blue)}
    \label{fig:remaining_sets}
\end{figure}

As a first summary we can conclude that the sequences repeat with the distance of \(2^{j}\) with \textit{j} being the number of even division by 2 steps. The set of natural numbers can be divided into subsets with individual solutions and cycles where not all of them can be solved immediately leaving room for the longer sequences as known from the problem (Figure \ref{fig:remaining_sets}).

Now it does make sense to rewrite the starting number \textit{n} as follows to track the remaining subsets in \(\mathbb{N}^{*}\) for which individual solutions need to be found:

\[n=2^{j}k+x; j,k,x \in \mathbb{N};n\in \mathbb{N}^{*}\]

where \textit{j} is the number of division steps (\textit{E}) in the sequence and \textit{x} is a specific number of the sequence and the lowest starting number showing the sequence. It is shown later that, in principle, an infinite sequence can be defined which would still lead to a final number lower than the initial starting number as well as a path towards infinity using alternating $OE$ steps.

Note that in the difference between starting and final value always grows with the value of $k$ (see Table \ref{tab:final_number}). This is of importance for the assessment of a potentially closed cycle with the same value for start and end of the sequence (see 1, 4, 2, 1). The subsequent starting values 5 and 9 have ending values of 4 and 7, respectively. This means that only the lowest number of any sequence can be the one which is a closed sequence. 

\subsection{Calculation of the first lower number in a sequence}
From the calculation of all steps in the sequence a formula for the first lowest number $y$ can be written in a general way as:

\[y=\frac{3^{i}\left ( 2^{j}k+x \right ) +m}{2^{j}}\]

With \textit{i} being the number of \textit{O} operations, \textit{j} the number of \textit{E} operations. \textit{x} is the starting number of the sequence and \textit{m} an adder following from the “+1” in the \textit{O} operation step formula.

\subsection{Subsequent lower number in the same cycle sequence}

Be \(y_{k}\) the first lower number for the sequence 

\[n_{k+1}=2^{j}(k+1)+x\]

is

\[y_{k+1}=y_{k}+3^{i}\]

\begin{table}
    \centering
    \begin{tabular}{|c|c|c|c|c|c|}
    \hline
 \multicolumn{5}{|c|}{Case \textit{j}=5, \textit{i} =3} & \\
 \hline
     \textit{k}=0    & \textit{k}=1 & \textit{k}=2 & \textit{k}=3 & \textit{k}=4 & Step no.\\
     \hline
 \multicolumn{5}{|c|}{\(2^{j}k+11= 32k+11\)  } & \\
 \hline
         11 & 43    & 75    & 107   & 139   & 1 \\
34 & 130   & 226   & 322   & 418   & 2 \\
17 & 65    & 113   & 161   & 209   & 3 \\
52 & 196   & 340   & 484   & 628   & 4 \\
26 & 98    & 170   & 242   & 314   & 5 \\
13 & 49    & 85    & 121   & 157   & 6 \\
40 & 148   & 256   & 364   & 472   & 7 \\
20 & 74    & 128   & 182   & 236   & 8 \\
10 & 37    & 64    & 91    & 118   &   \\
   & $10+27$ & $10+27\cdot2$ & $10+27\cdot3$ & $10+27\cdot4$ &   \\
   \hline
 \multicolumn{5}{|c|}{\(10+3^{i}k = 10+27k\)} & \\
         \hline
    \end{tabular}
    \caption{Sequences of the same type with a length of 8 steps}
    \label{tab:final_number}
\end{table}

\subsection{An example with the length 16}
Let’s assume we want to construct a sequence with the length 16. We know that this sequence has to be constructed with 10 \textit{E} operations and 6 \textit{O} operations to receive a lower number than the starting number with the 16th step.
The simplest (but not only) unique sequence with that length consist of alternating operations followed by E operations. Once this sequence has been identified it is possible to calculate the adder number \textit{m} for the general formula. Here the value is 665. The next step is to calculate the starting number offset \textit{x} which has to be a natural number and needs to fulfil the specific inequality for this sequence:

\[x> \frac{665}{2^{10}-3^{6} }=2.254\]

Further, the smallest natural number \textit{x} has to be found which gives a natural number result for the following stopping value diophantine equation:

\[y_{0}=    \frac{3^{6}x +665}{2^{10}}\]

or rewritten in the form of a diophantine equation as:

\[2^{10}y_{0} - {3^{6}x  - 665} = 0\]

or more general:

\[2^{E}y_{0} - {3^{O}x  - (\sum_{i=1}^{O-1}3^{O-1-i}2^{i-1})} = 0\]

where $O$ is the number of odd steps and $E$ the number of even steps. 

It is obvious that the equation has a solution as the requirement of a greatest common divisor (GCD) needs to be 1 is always fulfilled. The factors for $x$ and $y$ are powers of 2 and 3, respectively. The above equation gives the result \(x=575\) the final equation for the sequence is given by:

\[y=    \frac{3^{6}\left (  2^{10}k + 575 \right )+665}{2^{10}}\]

For $k=0$ the formula gives the result 410 a the first value lower than the starting value 575 which is in line with Table \ref{tab:my_label}.

\begin{table}
    \centering
    \begin{tabular}{|c|c|c|c|c|}
    \hline
       Step number  & Operation & Adder formula  & Adder value & Example \\
         \hline
         &  &  &  & 575\\
       1  &\textit{O}  &  \(3^{5}\cdot 2^{0}\)&243  & 1726\\
       2  & \textit{E} &  &  & 863\\
       3  & \textit{O} & \(3^{4}\cdot 2^{1}\) & 162 & 2590\\
       4  & \textit{E} &  &  & 1295\\
       5  & \textit{O} & \(3^{3}\cdot 2^{2}\) & 108 & 3886\\
       6  & \textit{E} &  &  & 1943\\
       7  & \textit{O} & \(3^{2}\cdot 2^{3}\) &  72& 5830\\
       8  &\textit{E}  &  &  & 2915\\
       9  & \textit{O} & \(3^{1}\cdot 2^{4}\) & 48 & 8746\\
       10  & \textit{E} &  &  & 4373\\
       11  & \textit{O} & \(3^{0}\cdot 2^{5}\) & 32 & 13120\\
        12 & \textit{E} &  &  & 6560\\
        13 & \textit{E} &  &  & 3580\\
        14 & \textit{E} &  &  & 1680\\
        15 & \textit{E} &  &  & 820\\
        16 & \textit{E} &  &  & 410\\
         &  &  & \textbf{665} & \\
         \hline
    \end{tabular}
    \caption{Base sequence with length 16 and derived example}
    \label{tab:my_label}
\end{table}

\begin{figure}
    \centering
    \includegraphics[width=0.9\linewidth]{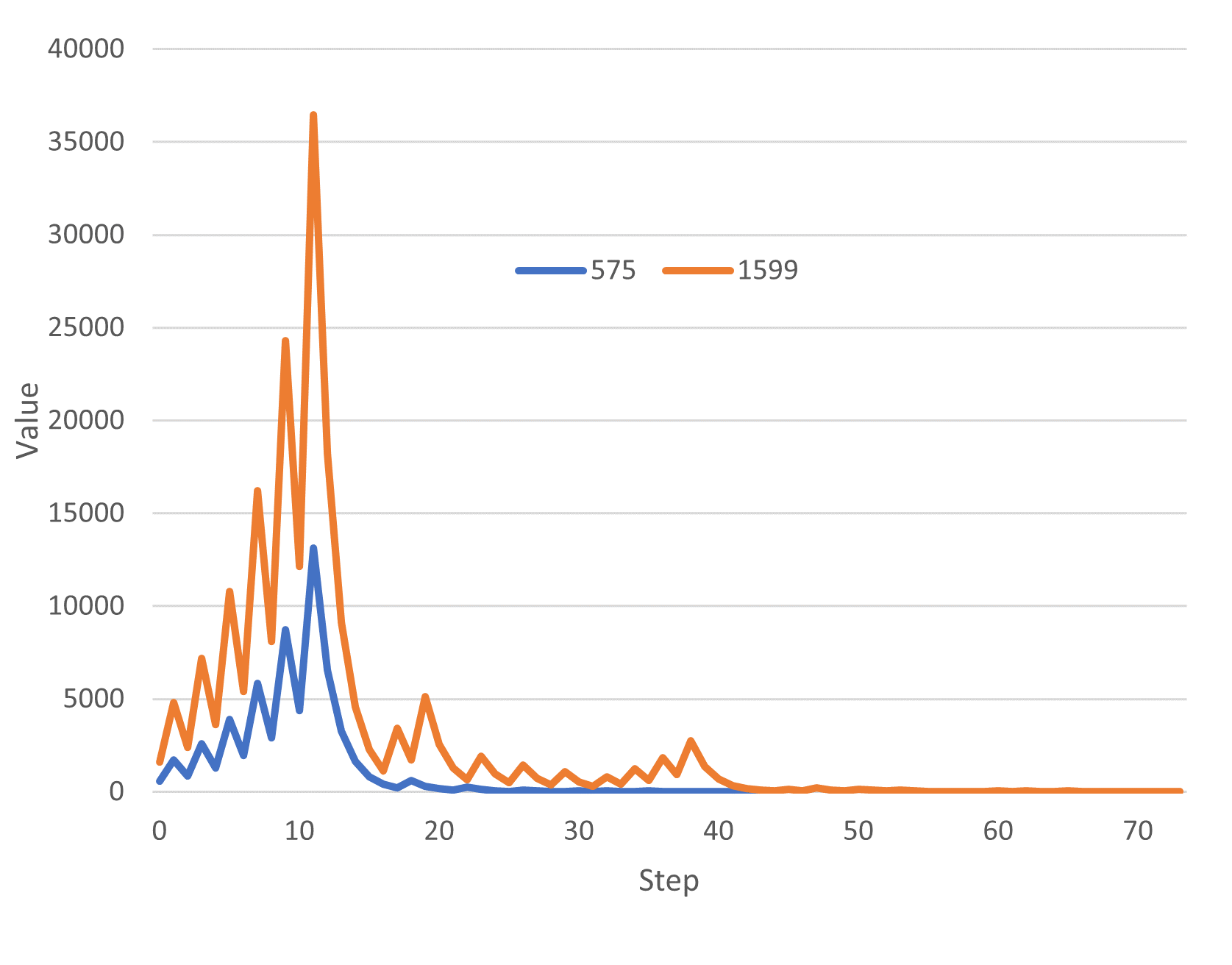}
    \caption{Collatz sequences for 575 and 1599}
    \label{fig:numbergraph}
\end{figure}

The full sequence until the final value 1 can be seen in Fig.\ref{fig:numbergraph}. Note the alternating steps in the beginning. The next number according to the equation (1599) shows exactly the same behaviour until the sequence ends and the new number defines the subsequent cycle.

\section{Only one closed cycle}

In an earlier section it has been shown that a general formula for a cycle with alternating $O$ and $E$ steps followed by a number of $E$ steps is given by:

\[2^{E}y_{0} - {3^{O}x  - (\sum_{i=1}^{O-1}3^{O-1-i}2^{i-1})} = 0\]

Assuming that a closed circle is finished with the same initial value $x$ this means that $y=x$ and the formula gets

\[2^{E}x - {3^{O}x  - (\sum_{i=1}^{O-1}3^{O-1-i}2^{i-1})} = 0\]

or
\[x = \frac{\sum_{i=1}^{O-1}3^{O-1-i}2^{i-1}}{2^{E} - 3^{O}}\]

The sum in the numerator can be rewritten as 

\[\sum_{i=1}^{O-1}3^{O-1-i}2^{i-1} = 3^{O} - 2^{O} \]

Hence, it can be concluded that the fraction needs to be an integer number and
\[ (2^{E} - 3^{O}) |(3^{O} - 2^{O})\]

The denominator $(3^{O} - 2^{O})$ is always positive for positive integers. To gain a positive result also the numerator $(2^{E} - 3^{O})$ needs to be positive, and $E > O \cdot log_2(3)$. However, as the power of $O$ is dominant in both numerator and denominator both powers $O$ and $E$ need to be in a close relation such as e.g. $E = O+1$ which is only the case for small values of $O$. One possible solution would be e.g.

\[ (O = 1) \wedge (E=2)  \]

where

\[x = \frac{3^{1} - 2^{1}}{2^{2} - 3^{1}} = \frac{3 - 2}{4 - 3} = 1\]

leading to the known sequence
\[ 1 \xrightarrow{O}4 \xrightarrow{E}2 \xrightarrow{E} 1\]

The above argument seems only true for the sequence of alternating $O$ and $E$ steps with subsequent $E$ steps as only this leads to the sum as well as the formula for the sum in the numerator. Other sequences might be addressed in the same manner but this is not entirely clear.

\section{27 and and more general cycles}
\subsection{The case of 27}

For the well known number 27 the sequence takes 96 steps to reach the next lower number - 23. The full sequence and the sequence for the subsequent number in this set can be found in the appendix. In short, the 96 steps consist of 59 \textit{E} operations and 37 \textit{O} operations leading to the following statements:
\begin{itemize}
    \item The distance between subsequent numbers with this sequence is  \(2^{59} = 576,460,752,303,423,488\)
    \item The next number with exactly the same sequence until its next lower number is 27 with the periodicity value added. The number is: 576,460,752,303,423,515
    \item The first smaller number of this subsequent number is \(3^{37} + 23 = 450,283,905,890,997,386 \)

\end{itemize}

\subsection{The "infinite" sequence}

It is clear now that, due to the nature of the sequence construction the set \(\mathbb{N}\) will only be filled with an infinite number of subsets of 2n and that always a longer sequence of \((OE)^{(j-times)}+E^{(i-j)-times}\) can be defined, as well as that the sequence according to the above formulas can be constructed.
The “infinite” sequence is of the following form:
\begin{itemize}
    \item \(((2^{i}-3^{j})> 0)\wedge ((2^{i}-3^{j})< 0)\wedge (i+j\rightarrow \infty )\)
    \item \textit{j} steps with \textit{OE} followed by (\textit{i}-\textit{j}) steps of \textit{E} 
    \item The adder constant for the first lower number formula is
    
\[m=\sum_{i=1}^{n}3^{n-i}2^{i-1}\]

\end{itemize}
Now let's consider the following sequence of alternating $O$ and $E$ steps and $n$ times similar to the previous sequences but without a number of $E$ steps to end up a number below the initial number.

With the number 575 we have again an adder value of 665 but this time the formula changes to:

\[2^{6}y_{0} - {3^{6}x  - 665} = 0\]

a general solution of this equation is:

\[y = 729t_{0} +174230 \]

\[x = 64t_{0} +15295\]

the smallest natural number fulfilling this equation is 63 which rises up to $y=1456$ at that point. The number from our example 575 rises up to $y = 13120$ according to the formula which is also in line with our observations. More precise, the formula result is $y=6560$ but it is clear that, following the definition of the sequence, the second to last number is double the value of $y$. The rational behind the cutoff is as follows. Consider an infinite sequence of $OE$ steps which ultimately leads to an infinite value of the final value. Will this one also have a solution of the equation and does this mean that this sequence will exist?

The next question would be how solutions of the diophantine equation behave with increasing sequence lengths. \ref{tab:starting_values} shows the result and that the initial lowest starting value depends on the sequence length $n$ and gives the result of $2^{n+1}$ which rises up to a maximum value within the sequence of $2 \cdot (3^{n}-1 )$.

\begin{table}
    \centering
    \begin{tabular}{|c|c|c|}
    \hline
       Steps of $OE$  & lowest starting number & highest value  \\
         \hline
       1  &1  &  4\\
       2   & 3 & 16 \\
       3  & 7 &  52\\
       4  & 15 & 160 \\
       5  &  31& 484 \\
       6  &  63& 1456 \\
       7  & 127 & 4372 \\
       8  & 255 & 13120 \\
       9  & 511 & 39364\\
       10  & 1023 & 118096 \\
$n$ & $2^{n+1}$ & $2 \cdot (3^{n}-1 )$\\
         \hline
    \end{tabular}
    \caption{Starting number and highest number after $n$ alternating $OE$ steps}
    \label{tab:starting_values}
\end{table}

\begin{figure}
    \centering
    \includegraphics[width=0.9\linewidth]{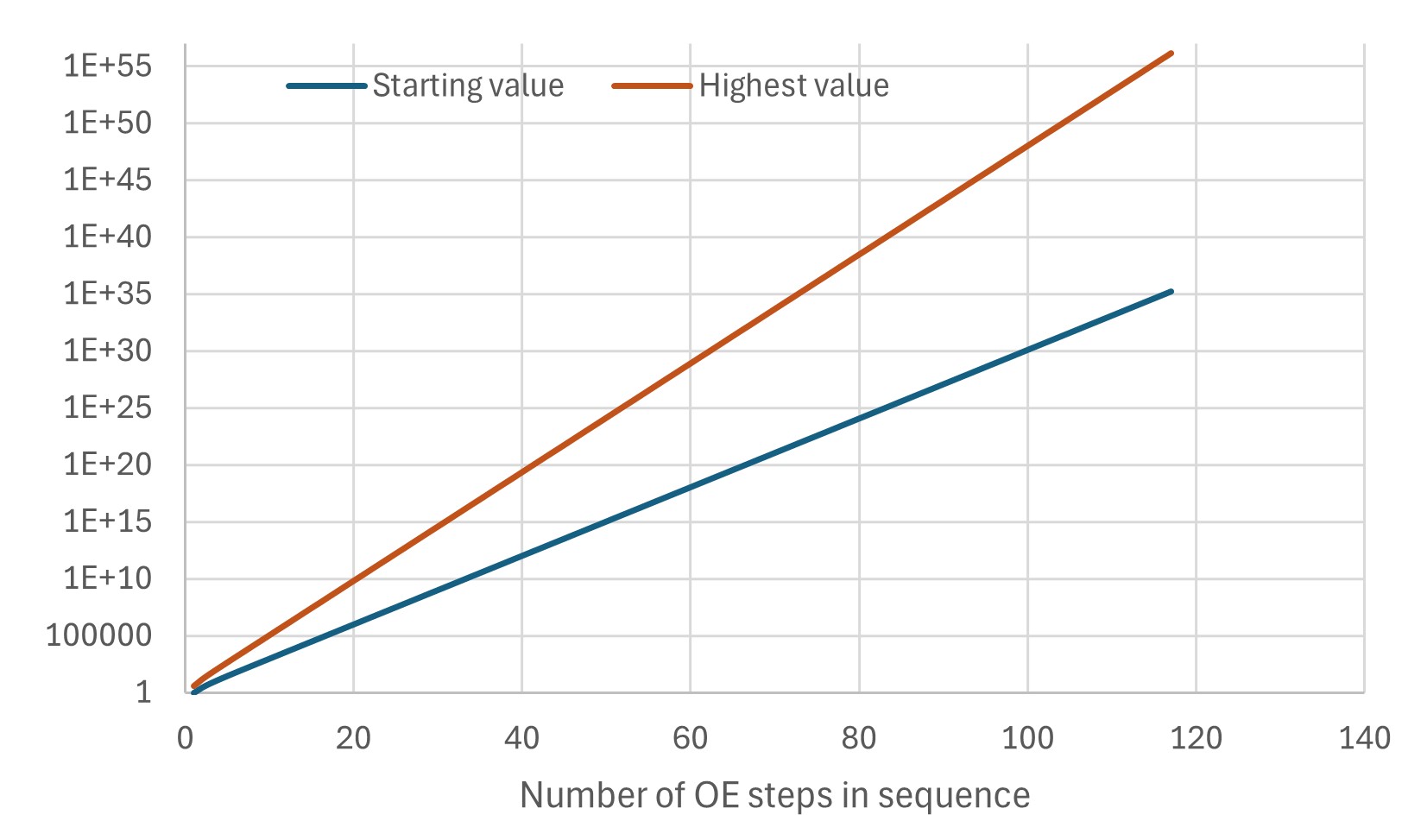}
    \caption{Lowest number to fulfill the equation for each $OE$ steps length and respective highest value }
    \label{fig:rising_numbers}
\end{figure}

\ref{fig:rising_numbers} shows the increase of starting value and the respective highest value during this sequence on a logarithmic scale. Certainly, as discussed earlier several statements are true:

\begin{itemize}
    \item For each sequence length a diophantine equation can be formed which will also have a solution as the the GCD of both factors is always 1
    \item A simple repetitive $OE$ sequence will also have a solution and it could be shown that there is also a solution to be expected when the sequence length goes towards infinity according to the respective formula
    \item It is also true that also the starting value will go towards infinity with $2^{n+1}$
\end{itemize}

An interesting question connected to the periodicity of starting values $2^ik$ is as follows:
If we consider the statement that all patterns which are the same have this periodicity and that otherwise there needs to be a difference in the sequence of steps. Would an infinite $OE$ sequence be periodic with the subsequent patterns of $3n+1$ and would this be a contradiction as the period should also become infinite? And if this is the case, could this be true for all sequences?

\subsection{The solved subsets of $2^ {E}+j$}

\begin{figure}
    \centering
    \includegraphics[width=0.9\linewidth]{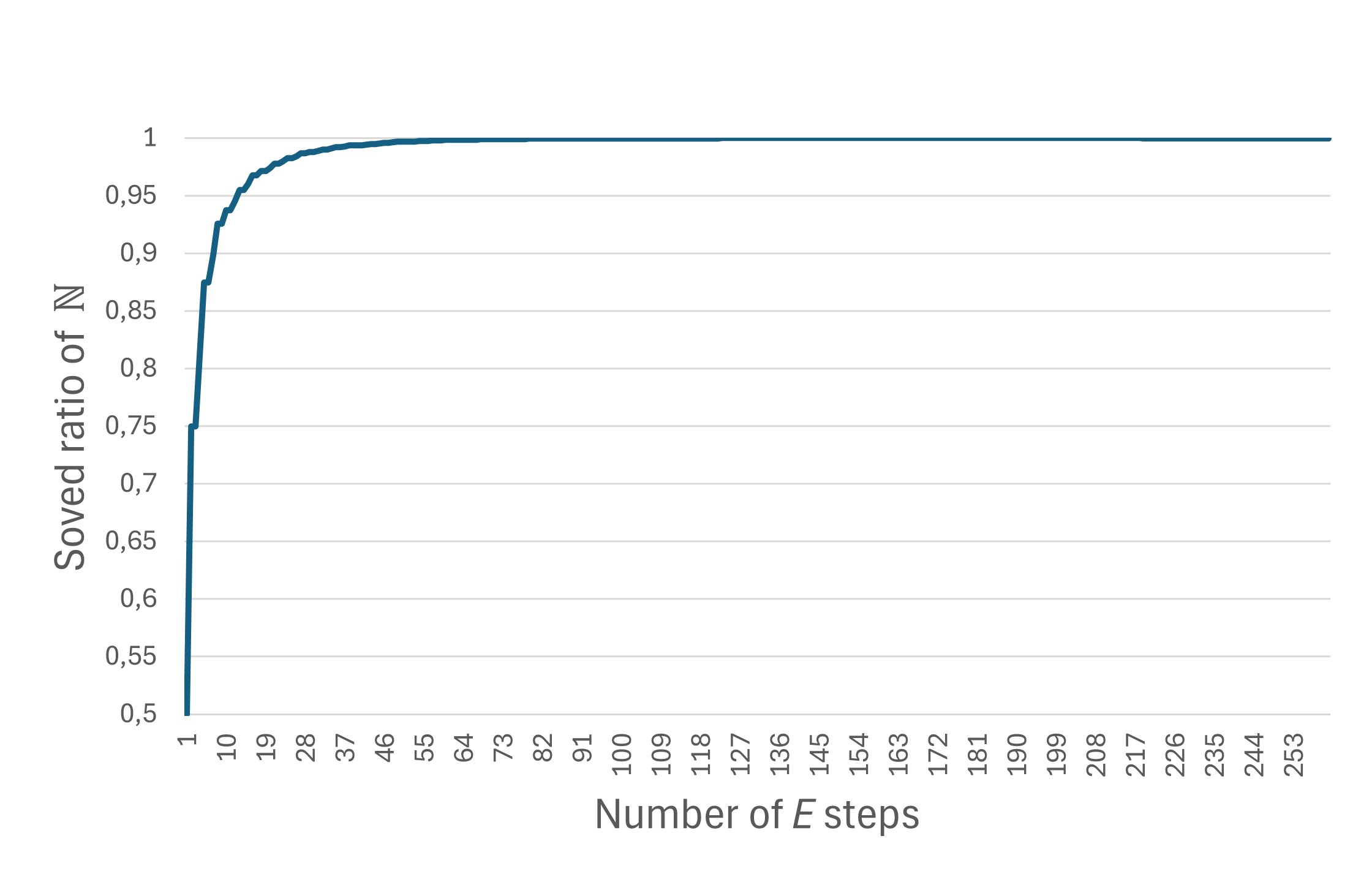}
    \caption{Ratio of natural numbers leading to a first smaller number latest after E steps}
    \label{fig:ratio_of_N}
\end{figure}

As shown earlier all numbers within the same subset of $2^ {E}+j$ show the same behavior especially in the way they end up at the first smaller number than the initial value. One could conclude that e.g. after the first step all even numbers are solved as already the first value is $n/2$, while all odd numbers remain unsolved as they cannot lead to a smaller number after one $E$ step. This means that a probability of $0.5$ of all integers are solved after a sequence with 1 $E$ step and $0.5$ of the values remain in an unsolved subset which is doubled in the next step.
Here, $4k+1$ and $4k+3$ remain where the first subset is solved after 2 $E$ steps ($OEE$) and the latter is not. The solving rate is therefore increased to $0.75$ with one unsolved subset remaining. It has been shown already that there is no possible solution with 3 $E$ steps, however, the number of subsets needs to be doubled for the next step as the sequence of $OE$ is different. This means 2 subsets are unsolved leading to 4 subsets to be investigated with 4 $E$ steps. Here, only one subset $2^4k+3$ is solved with $0.8125$ of the natural numbers solved and 3 subsets remain unsolved leading to 6 more unsolved subsets in the next step. It is therefore interesting to see if the number of unsolved subsets grows faster than the number of solved subsets leading to a saturation at at number smaller than 1. Figure \ref{fig:ratio_of_N} shows this behavior and indeed that the the ratio approaches values close to unity.

\section{Conclusion}
It could be shown that each starting number in the Collatz problem can be rewritten in the form \(2^{j}k+x\) and a specific sequence of steps can be attributed following the Collatz function definition leading to the first value lower than the starting number. It further could be shown that, given a specific length of such a sequence and provided that this length can lead to a lower number at all at least one sequence can be defined and the respective parameter values can be calculated. The subsequent sequence after this first one will most likely not be of the same type again, which can also be seen easily, leading to the well observed differences in total stopping time.
With respect to the proof of the Collatz problem it seems plausible to follow the idea that it is sufficient to show that all numbers (exept 1) will lead to a lower number and thus need to end up with 1. The sequences are identified needed to fill the set of natural numbers and identify sets based on their offset value \textit{x} as this is due to \(2^{j}0+x=x\) in any case the first starting value. Afterwards all numbers which are part of this set are solved immediately.
In terms of visualization the subset tree (see Figure \ref{fig:remaining_sets}) will also be of interest considering a larger depth of $2^{j}$.

It further could be shown that an infinite sequence of $OE$ steps will have a solution and will also go towards infinity in both initial starting value and highest value in this sequence. It needs to be discussed if this idea can serve towards a proof. Following the number of solved and newly generated subsets with each step it seems that the solutions will ultimately lead to a saturation at a probability at $1$ which could be interpreted as there is a zero probability that infinite Collatz sequences exists. 

Future work could focus on the following key points:
\begin{itemize}
    \item Identify an algorithm to be able to attribute a sequence length to a given starting number without running through the Collatz formula
    \item How can the number \textit{n} be excluded from further investigation fast in case the number is already part of an identified sequence.
    \item Can it be proven that the probability of subsets been solved will approach unity.
    \item Is a periodic behavior in contradiction with infinite sequence lengths.
    
\end{itemize}


%

\newpage
\appendices
\section{Examples for the 2 sequences with length 8}

\begin{table}[H]
    \centering
    \begin{tabular}{|c|c|c|c|c|}
    \hline
Number 1 & Number 2 & Number 3 & Step & Cycle \\
\hline
11       & 43       & 331      & 1    & O     \\
34       & 130      & 994      & 2    & E     \\
17       & 65       & 497      & 3    & O     \\
52       & 196      & 1492     & 4    & E     \\
26       & 98       & 746      & 5    & E     \\
13       & 49       & 373      & 6    & O     \\
40       & 148      & 1120     & 7    & E     \\
20       & 74       & 560      & 8    & E     \\
10       & 37       & 280      &      &    \\
\hline
    \end{tabular}
    \caption{Sequence \(2^{5}k+11\)}
    \label{tab:cycle_8_1}
\end{table}

\begin{table}[H]
    \centering
    \begin{tabular}{|c|c|c|c|c|}
    \hline
Number 1 & Number 2 & Number 3 & Step & Cycle \\
\hline
23       & 55       & 343      & 1    & O     \\
70       & 166      & 1030     & 2    & E     \\
35       & 83       & 515      & 3    & O     \\
106      & 250      & 1546     & 4    & E     \\
53       & 125      & 773      & 5    & O     \\
160      & 376      & 2320     & 6    & E     \\
80       & 188      & 1160     & 7    & E     \\
40       & 94       & 580      & 8    & E     \\
20       & 47       & 290      &      &      \\
\hline
    \end{tabular}
    \caption{Sequence \(2^{5}k+23\)}
    \label{tab:cycle_8_2}
\end{table}

\newpage

\section{The table with the sequence 27}

\begin{table}[H]
    \centering
    \begin{tabular}{|c|c|c|c|c|}
    \hline
Step & Value & Ops & Adder value & Subsequent \\
\hline
1        & 27       & \textit{O}        & 1,50095E+17 & 576460752303423515  \\
2        & 82       & \textit{E}        &             & 1729382256910270546 \\
3        & 41       & \textit{O}        & 1,00063E+17 & 864691128455135273  \\
4        & 124      & \textit{E}        &             & 2594073385365405820 \\
5        & 62       & \textit{E}        &             & 1297036692682702910 \\
6        & 31       & \textit{O}        & 1,33417E+17 & 648518346341351455  \\
7        & 94       & \textit{E}        &             & 1945555039024054366 \\
8        & 47       & \textit{O}        & 8,8945E+16  & 972777519512027183  \\
9        & 142      & \textit{E}        &             & 2918332558536081550 \\
10       & 71       & \textit{O}        & 5,92966E+16 & 1459166279268040775 \\
11       & 214      & \textit{E}        &             & 4377498837804122326 \\
12       & 107      & \textit{O}        & 3,95311E+16 & 2188749418902061163 \\
13       & 322      & \textit{E}        &             & 6566248256706183490 \\
14       & 161      & \textit{O}        & 2,63541E+16 & 3283124128353091745 \\
15       & 484      & \textit{E}        &             & 9849372385059275236 \\
16       & 242      & \textit{E}        &             & 4924686192529637618 \\
17       & 121      & \textit{O}        & 3,51388E+16 & 2462343096264818809 \\
18       & 364      & \textit{E}        &             & 7387029288794456428\\
19       & 182      & \textit{E}        &             & 3693514644397228214 \\
20       & 91       & \textit{O}        & 4,68517E+16 & 1846757322198614107 \\
21       & 274      & \textit{E}        &             & 5540271966595842322 \\
22       & 137      & \textit{O}        & 3,12344E+16 & 2770135983297921161 \\
23       & 412      & \textit{E}        &             & 8310407949893763484 \\
24       & 206      & \textit{E}        &             & 4155203974946881742 \\
25       & 103      & \textit{O}        & 4,16459E+16 & 2077601987473440871 \\
26       & 310      & \textit{E}        &             & 6232805962420322614 \\
27       & 155      & \textit{O}        & 2,7764E+16  & 3116402981210161307\\
28       & 466      & \textit{E}        &             & 9349208943630483922  \\
29       & 233      & \textit{O}        & 1,85093E+16 & 4674604471815241961  \\
30       & 700      & \textit{E}        &             & 14023813415445725884 \\
31       & 350      & \textit{E}        &             & 7011906707722862942  \\
32       & 175      & \textit{O}        & 2,46791E+16 & 3505953353861431471  \\
33       & 526      & \textit{E}        &             & 10517860061584294414 \\
34       & 263      & \textit{O}        & 1,64527E+16 & 5258930030792147207  \\
35       & 790      & \textit{E}        &             & 15776790092376441622 \\
36       & 395      & \textit{O}        & 1,09685E+16 & 7888395046188220811  \\
37       & 1186     & \textit{E}        &             & 23665185138564662434 \\
38       & 593      & \textit{O}        & 7,31232E+15 & 11832592569282331217 \\
39       & 1780     & \textit{E}        &             & 35497777707846993652 \\
40       & 890      & \textit{E}        &             & 17748888853923496826 \\
41       & 445      & \textit{O}        & 9,74976E+15 & 8874444426961748413  \\
42       & 1336     & \textit{E}        &             & 26623333280885245240 \\
43       & 668      & \textit{E}        &             & 13311666640442622620 \\
44       & 334      & \textit{E}        &             & 6655833320221311310 \\
45       & 167      & \textit{O}        & 2,59993E+16 & 3327916660110655655  \\
46       & 502      & \textit{E}        &             & 9983749980331966966  \\
47       & 251      & \textit{O}        & 1,73329E+16 & 4991874990165983483  \\
48       & 754      & \textit{E}        &             & 14975624970497950450 \\
49       & 377      & \textit{O}        & 1,15553E+16 & 7487812485248975225  \\
50       & 1132     & \textit{E}        &             & 22463437455746925676 \\

\hline
    \end{tabular}
    \caption{Sequence for 27 and subsequent number)}
    \label{tab:cycle_27}
\end{table}

\begin{table}
    \centering
    \begin{tabular}{|c|c|c|c|c|}
    \hline
Step & Value & Ops & Adder value & Subsequent \\
\hline
51       & 566      & \textit{E}       &             & 11231718727873462838 \\
52       & 283      & \textit{O}        & 1,5407E+16  & 5615859363936731419  \\
53       & 850      & \textit{E}        &             & 16847578091810194258 \\
54       & 425      & \textit{O}        & 1,02713E+16 & 8423789045905097129  \\
55       & 1276     & \textit{E}        &             & 25271367137715291388 \\
56       & 638      & \textit{E}        &             & 12635683568857645694 \\
57       & 319      & \textit{O}        & 1,36951E+16 & 6317841784428822847  \\
58       & 958      & \textit{E}        &             & 18953525353286468542 \\
59       & 479      & \textit{O}       & 9,13009E+15 & 9476762676643234271  \\
60       & 1438     & \textit{E}        &             & 28430288029929702814 \\
61       & 719      & \textit{O}       & 6,08672E+15 & 14215144014964851407 \\
62       & 2158     & \textit{E}        &             & 42645432044894554222 \\
63       & 1079     & \textit{O}        & 4,05782E+15 & 21322716022447277111 \\
64       & 3238     & \textit{E}      &             & 63968148067341831334 \\
65       & 1619     & \textit{O}       & 2,70521E+15 & 31984074033670915667 \\
66       & 4858     & \textit{E}        &             & 95952222101012747002 \\
67       & 2429     & \textit{O}      & 1,80347E+15 & 47976111050506373501 \\
68       & 7288     & \textit{E}        &             & 143928333151519120504 \\
69       & 3644     & \textit{E}        &             & 71964166575759560252  \\
70       & 1822     & \textit{E}        &             & 35982083287879780126  \\
71       & 911      & \textit{O}        & 4,80926E+15 & 17991041643939890063  \\
72       & 2734     & \textit{E}       &             & 53973124931819670190  \\
73       & 1367     & \textit{O}       & 3,20618E+15 & 26986562465909835095  \\
74       & 4102     & \textit{E}        &             & 80959687397729505286  \\
75       & 2051     & \textit{O}        & 2,13745E+15 & 40479843698864752643  \\
76       & 6154     & \textit{E}        &             & 121439531096594257930 \\
78       & 9232     & \textit{E}       &             & 182159296644891386896 \\
79       & 4616     & \textit{E}       &             & 91079648322445693448  \\
80       & 2308     & \textit{E}        &             & 45539824161222846724  \\
81       & 1154     & \textit{E}        &             & 22769912080611423362  \\
82       & 577      &\textit{ O}        & 7,59982E+15 & 11384956040305711681  \\
83       & 1732     & \textit{E}        &             & 34154868120917135044  \\
84       & 866      & \textit{E}        &             & 17077434060458567522  \\
85       & 433      & \textit{O}        & 1,01331E+16 & 8538717030229283761  \\
86       & 1300     & \textit{E}       &             & 25616151090687851284 \\
87       & 650      & \textit{E}       &             & 12808075545343925642 \\
88       & 325      & \textit{O}        & 1,35108E+16 & 6404037772671962821  \\
89       & 976      & \textit{E}        &             & 19212113318015888464 \\
90       & 488      & \textit{E}        &             & 9606056659007944232  \\
91       & 244      & \textit{E}        &             & 4803028329503972116  \\
92       & 122      & \textit{E}        &             & 2401514164751986058  \\
93       & 61       & \textit{O}        & 7,20576E+16 & 1200757082375993029  \\
94       & 184      & \textit{E}        &             & 3602271247127979088 \\
95       & 92       & \textit{E}        &             & 1801135623563989544  \\
96       & 46       & \textit{E}        &             & 900567811781994772   \\
          & 23       &           &             & \textbf{450283905890997386}  \\
          \hline
Total & 96       & Adder & 1,10093E+18 &                     \\
O steps   & 37       &           &             &                     \\
E steps   & 59       &           &             &                     \\
\hline
    \end{tabular}
    \caption{Sequence for 27 and subsequent number - continued)}
    \label{tab:cycle_27b}
\end{table}




\ifCLASSOPTIONcaptionsoff
  \newpage
\fi



%
\newpage
\bibliography{bibtex/bib/Collatz}
\bibliographystyle{ieeetr}

to be added


%








\end{document}